\documentclass[journal]{IEEEtran}
\usepackage{amsmath,amsthm}
\usepackage{amsfonts,amssymb}
\usepackage{CJK}
\usepackage{multicol}
\usepackage{multirow}
\usepackage{diagbox}
\usepackage{graphicx}
\usepackage{subfig} 
\usepackage{lineno} 
\allowdisplaybreaks 
\usepackage{epstopdf} 
\usepackage{pgf,tikz}
\usepackage{tikz}
\usepackage{graphicx}     
\usepackage{subfig} 
\usepackage{subfloat}
\usepackage{cite}
\usepackage[subfigure]{tocloft}
\usetikzlibrary{arrows,shapes,snakes,shadows,positioning,automata,patterns}
\usetikzlibrary{trees,decorations.pathmorphing,decorations.markings}
\usetikzlibrary{backgrounds}
\usepackage[
            pdfstartview=FitH,
            CJKbookmarks=true,
            bookmarksnumbered=true,
            bookmarksopen=true,
            colorlinks, 
            pdfborder=001,   
            linkcolor=blue,
            anchorcolor=blue,
            citecolor=blue
            ]{hyperref}
\newtheorem{theorem}{Theorem}
\newtheorem{assumption}{Assumption}

\newtheorem{lemma}{Lemma}

\newtheorem{remark}{Remark}

\ifCLASSINFOpdf
\else
\fi

\begin{document}

\title{Novel Distributed Algorithms Design for Nonsmooth Resource Allocation  on Weight-Balanced Digraphs}

\author{Xiaohong Nian, Fan Li,  Dongxin Liu
\thanks{This work is supported in part by the National Natural Science Foundation
	of China under Grant 62173347. (Corresponding author: Fan Li)}
\thanks{X. Nian, F. Li, and D. Liu are with Key Lab of Institute of cluster unmanned
	system, School of automation, Central South University, Changsha, 410083, China  (e-mail:funnynice@csu.edu.cn).}
}
\maketitle

\begin{abstract}
In this paper, the distributed resource allocation problem on strongly connected and weight-balanced digraphs is investigated, where the decisions of each agent are restricted to satisfy the coupled network resource constraints  and heterogeneous general convex sets.  
Moreover, the local cost function can be non-smooth. In order to achieve the exact optimum of the nonsmooth resource allocation problem, a novel continuous-time distributed algorithm based on the gradient descent scheme and differentiated projection operators is proposed. 
With the help of the set-valued LaSalle invariance principle and nonsmooth analysis, it is demonstrated that the algorithm converges asymptotically to the global optimal allocation. Moreover, for the situation where local constraints are not involved and the cost functions are differentiable with Lipschitz gradients, the convergence of the algorithm to the exact optimal solution  is exponentially fast. Finally, the effectiveness of the proposed algorithms is illustrated by simulation examples.
\end{abstract}


\IEEEpeerreviewmaketitle

\section{Introduction}
Distributed resource allocation is an essential issue in distributed optimization problems with widespread applications in many domains, such as smart grids \cite{grad1,grid2}, wireless networks \cite{net1,net2}, transportation \cite{lu1,lu2}, and multi-robot networks \cite{rob1}. As a process of allocating resources among agents, the core objective of the resource allocation problem is to seek the optimal decision (or optimal allocation) for each agent to minimize the sum of the costs of all agents under the limited resource constraints.

The distributed approach (see \cite{res1,res2,LinwenTing2021,Mag}) to solving the resource allocation problem has received increasing attention owing to its ease of scaling in large-scale networks as well as better robustness and privacy preservation than the conventional centralized approach. 
Up to now, despite discrete-time algorithms have been studied in most of the literature (see \cite{Lisan,Lisan2,Lisan3,Xu,Ding,ZhangTsp}), owing to the fact that continuous-time algorithms are more adaptable to be implemented in physical systems, and more convenient to be investigated for convergence in terms of analytical tools in control theory, distributed continuous-time resource allocation algorithms have attracted significant attention as well (see    \cite{Yi2016,Xuefang2019,Wangex2018,Chen2021TSM,2018Dis,YinAuto2019,Lgriad,He2017,Sunlanlan2021,KiaAuto2021}).
In \cite{Yi2016}, the primal-dual projection-based algorithms were designed to resolve the smooth resource allocation problem on the undirected communication topology. The application of the algorithms to the economic dispatch problem in power systems was also further discussed. Later, in \cite{Xuefang2019}, the distributed resource allocation problem with alternating  ``small" and ``large" communication delays is studied, and a switching algorithm is proposed using switching techniques. Moreover, the Lyapunov functional theory is employed to determine the sufficient conditions for the exponential convergence of the algorithm related to the time delay.
Besides, in \cite{Wangex2018}, an extremum-seeking algorithm first-order dynamics is developed and extended to a second-order algorithm equipped with low-pass filters for resolving the resource allocation problem under undirected graphs, and the algorithm can converge within the neighborhood of optimal solution. In \cite{Chen2021TSM}, the method of $\epsilon$-exact penalty is adopted to deal with the box-like constraint, and a class of distributed optimization algorithms with predefined time convergence is presented based on the TBG scheme. Moreover, in \cite{2018Dis}, the resource allocation problem with box constraint is handled by means of the $\theta$-logarithmic barrier penalty functions. In \cite{YinAuto2019}, an initialization-free distributed algorithm is proposed based on the dual gradient method, where the agent interacts with its neighbors only for dual variables information, thereby effectively avoiding privacy leakage.
In \cite{Lgriad,He2017}, the distributed Laplacian-gradient dynamics are developed coupled with nonsmooth exact penalty functions to deal with economic dispatch problems of first-order and second-order systems, respectively. In \cite{Sunlanlan2021}, the event-triggered scheme and the periodic communication scheme are introduced to reduce the communication burden with a discrete communication mechanism. In \cite{KiaAuto2021}, the algorithm is presented by virtue of the augmented Lagrangian function, and the same method as in \cite{Chen2021TSM} is adopted to deal with the box set constraint. It is worth pointing out that the above-mentioned algorithms in the literature      \cite{Yi2016,Xuefang2019,Wangex2018,Chen2021TSM,2018Dis,YinAuto2019,Lgriad,He2017,Sunlanlan2021,KiaAuto2021} all ask the cost functions to be differentiable, which is hardly realized in practical engineering problems.

In fact, in practical applications, the optimization problems with nonsmooth cost functions exist widely, such as the bandwidth allocation problem for wireless networks where the utility function of the service quality is nonsmooth (see \cite{Unon}). While most of the known algorithms cannot be applied to the non-smooth resource allocation problem for the restriction on the differentiability of cost functions. Recently, some literature has studied nonsmooth resource allocation problems with a subgradient-based approach (see \cite{non108,non109,Lian22021,Zeng2016}). In \cite{non108}, a distributed iterative strategy is designed based on a centerless algorithm, although it will fail when the agents are subject to local feasible constraints.
In \cite{non109}, a distributed Lagrangian algorithm is presented to deal with the resource allocation problem on an undirected graph. In \cite{Lian22021}, a distributed algorithm with an adaptive gain is designed with differential inclusion and the Laplacian-gradient method. By adopting a distance-based exact penalty function to handle the local feasible set constraint, the algorithm can tackle the nonsmooth resource allocation problem on undirected graphs. Besides, in \cite{Zeng2016}, the projection-based algorithm is developed via a gradient descent scheme to solve the nonsmooth resource allocation problem.

It should be mentioned that the aforementioned literature \cite{Yi2016,Xuefang2019,Chen2021TSM,non108,non109,Lian22021,Zeng2016} all required communication among agents through an undirected graph. However, due to practical environmental concerns, information transmission cannot be completely symmetric, and pursuing such bidirectional communication will incur high-cost expenses. In contrast, weight-balanced digraphs, as a generalized form of undirected graphs, are easier to realize in practical applications. Therefore, the problem of resource allocation on weight-balanced digraphs has an important research significance.

Motivated by the above observations, this paper is devoted to tackling the nonsmooth resource allocation problem with general set constraints. Moreover, despite that the same problem considered in this paper as in \cite{Deng2018,Zhu2019TSM,Deng2020}, we design completely different novel distributed algorithms. Compared to the existing results, our main contributions are as follows.
\begin{enumerate}
	\item We investigate the distributed resource allocation problem on weight-balanced digraphs, where all agents are subject to heterogeneous local constraints depicted by general convex sets and coupled network resource constraints. Moreover, the local cost function is allowed to be non-smooth.
	Even though similar problems have been studied in \cite{Yi2016,Xuefang2019,Wangex2018,Chen2021TSM,2018Dis,YinAuto2019,Lgriad,He2017,Sunlanlan2021,KiaAuto2021,non108,non109,Lian22021,Zeng2016}, a variety of drawbacks distinguish them from this paper. For example, in \cite{Yi2016,Xuefang2019,Wangex2018,Chen2021TSM,2018Dis,YinAuto2019,Lgriad,He2017,Sunlanlan2021,KiaAuto2021}, the cost function needs to be differentiable or even twice differentiable.  In \cite{2018Dis,Chen2021TSM,Lgriad,He2017}, there is no local feasible set considered, or just box constraint instead of general convex set constraint. In \cite{Yi2016,Xuefang2019,Chen2021TSM,non108,non109,Lian22021,Zeng2016}, only the case of undirected graphs which are the special instance of weight-balanced digraphs is studied.
	
	\item The gradient descent scheme and differential inclusion theory are adopted to design a navel distributed algorithm different from \cite{Yi2016,Deng2018,Zhu2019TSM}, in which an overall information packaged in the form of a sum of several local information is sent between agents. Furthermore, the differential projection operator is employed to address convex set constraints, which will enable the states of agents to always be within the set constraints. The convergence of the algorithm is analyzed by the set-valued LaSalle invariance principle.

	\item When encountering distributed resource allocation problems without local set constraints, under the assumption that the local cost function has a Lipschitz gradient, the algorithm has exponential convergence to the optimal  solution of the problem globally on strongly connected and weight-balanced digraphs.
	
\end{enumerate}

We organize the remaining sections of this paper as follows. 
In Section \ref{sec.basis}, some fundamental preliminaries are provided and the problem formulation is introduced. 
In Section  \ref{sec.r}, the main results are presented, in which two distributed algorithms are designed for resolving resource allocation problems and the convergence is rigorously proved. 
In Section \ref{sec.s}, some numerical simulations are given to substantiate the validity of algorithms. 
Finally, in Section  \ref{sec.c}, concluding remarks are summarized.

\textsl{Notations:}
$\mathbb{R}^n$ represents the $n$-dimensional real  Euclidean space. 
For a vector $x \in \mathbb{R}^n$, $\|x\|$ is the Euclidean norm and $x_i$ is the $i$th element. $col(x_1, \ldots, x_N)=[x_1^T, \dots ,x_N^T]^T.$ $\|A\|$ indicates the spectral norm of matrix $A$. $\otimes$ represents the kronecker product. $0_n(1_n)$ utilized the $n\times1$ zeros(ones) vector. $I_n$ is used to denote the $n \times n$ identity matrix. 
$B(x,\varsigma)=\{y \in \mathbb{R}^n ~|~ \|y-x\| < \varsigma\}.$
$\partial \Omega$ and $int(\Omega) $ are taken to represent the boundary set and the relative interiors of the set $\Omega \subset \mathbb{R}^n.$

\section{Preliminaries and Formulation}\label{sec.basis}

\subsection{Graph Theory}
Consider a directed graph  $\mathcal{G}=\left( \mathcal{V}, \mathcal{E},\mathcal{A} \right) $
utilized to represent the interaction topology of a $N$ agents network. $\mathcal{G}=\left( \mathcal{V}, \mathcal{E},\mathcal{A} \right) $ is specified by three parameters: the node set $\mathcal{V}=\{1,\dots ,N\}$ , the edge set  $
\mathcal{E}\in \mathcal{V}\times \mathcal{V}$ and the weighted adjacency matrix $\mathcal{A}=\left( a_{ij} \right) _{N\times N}\in \mathbb{R}^{N\times N}$. Agent $i$ can obtain information from agent $j$ while $e_{ij}\in  \mathcal{E}$. If $e_{ij}\in  \mathcal{E}$, then $a_{ij}>0$; otherwise, $a_{ij}=0 \left( a_{ii}=0, \forall \;  i \in \mathcal{V} \right)$. The directed path from agent $i$ to agent $j$ is defined by a sequence  of edges $\left( i,i_1 \right) \left( i_1,i_2 \right) ,\cdots ,\left( i_k,j \right) $.  The digraph is said to be strongly connected if any two agents are connected by at least one directed path.  Besides,  $d_{in}^{i}=\sum_{j=1}^N{a_{ij}}$  and $d_{out}^{i}=\sum_{j=1}^N{a_{ji}}$ are denoted to represent the in-degree and  out-degree of agent $i$, respctively. Correspondingly, the Laplacian matrix associate with $\mathcal{G}$ is formulated as $L= D^{in}-\mathcal{A}$ with $D_{in}=\text{diag}\left\{ d_{in}^{1},\dots ,d_{in}^{N} \right\} \in \mathbb{R}^{N\times N}$.
In particular, $\mathcal{G}$ is known as a weight-balanced digraph if $d_{in}^i=d_{out}^i $ for all $i \in \mathcal{V}$.

By denoting $Sym(L)=\frac{L+L^T}{2 \,}$, some equivalent descriptions about weight-balanced digraphs are expressed as follows:
\begin{itemize}
	\item $1_N^T L =0_N^T$;
	\item $Sym(L)$ is positive semidefinite;
	\item $\mathcal{G}$ is weight-balanced.
\end{itemize}

Use ${\hat{\lambda} _1}, \dots, \hat{\lambda} _N$  with $\hat{\lambda} _i \le \hat{\lambda} _j $ for $i < j$ to represent the eigenvalues of $Sym(L)$.  For a digraph $\mathcal{G}$, the strong connectivity implies that 0 is a simple eigenvalue of $Sym(L)$ and the other eigenvalues have positive real parts. For detailed introduction, it can be found in \cite{Godsil01,Rikos14}.

\begin{lemma}\label{lem.tu} (see \cite{Yang2019TAC})
	If the graph $\mathcal{G}$ is weight-balanced, then the Laplacian matrix $L \in \mathbb{R}^{n\times n}$ can be factored as $L=T diag\{0,J\}T^T$, where $J \in \mathbb{R}^{(N-1)\times(N-1)}$, by an orthogonal matrix
	$T=[r \; \;   R] \in \mathbb{R}^{N\times N}$ with 
	$r=\frac{1}{\sqrt{N}}1_N$. In particular, there hold   
	$r^TR=0_{N}^{T},\ R^TR=I_{N-1}$, and 
	$RR^T=I_N-\frac{1}{N}1_N1_{N}^{T}$.
\end{lemma}

\subsection{Convex Analysis and Projection}
The following definitions are collected from \cite{Rockafellar70}.

A function $f: \Omega  \to \mathbb{R}$ is called to be  convex if  its domain $\Omega \subset \mathbb{R}^n$  is convex and  
$f(\alpha u + (1-\alpha)v) \le \alpha f(u) + (1-\alpha)f(v), \,\forall  u, v \in  \Omega$ with $0 \le \alpha \le 1$. 
The subgradient of  $f : \mathbb{R}^n \to \mathbb{R}$  at $u$ is defined by the vector $\gamma \in \mathbb{R}^n$
satisfying $f(v) \geq f(u) +\gamma^T (v-u), \forall v \in \mathbb{R}^n$.
Furthermore, the subdifferential of $f$ at $u$ composed of all subgradient is denoted as $\partial f(u) $.
A function $f: \Omega  \to \mathbb{R}$ is $\omega$-strongly convex with $\omega >0$ if $(u-v)^T(\gamma_1-\gamma_2) \ge  \omega {\left\| {u - v} \right\|^2}, \forall \, u,v \in \Omega, \gamma_1 \in \partial f(u), \gamma_1 \in \partial f(v)$.
A function $f: \mathbb{R}^n  \to \mathbb{R}$ is said to be Lipschitz with constant $\theta$ if $
\left\| f\left( u \right) -f\left( v \right) \right\|\le \theta \left\| u-v \right\|, \,  \forall \, u,v \in \mathbb{R}^n$.

The normal cone to  a  closed convex set $\Omega$ at $u\in \Omega$  is given by $
\mathcal{N}_{\varOmega}\left( u \right) =\left\{ v\in \mathbb{R}^n:v^T\left( w-u \right) \le 0,\ \forall w\in \varOmega \right\} $.
The polar cone of $\mathcal{N}_{\varOmega}$ defined by $T_\Omega(u)=\{v \in \mathbb{R}^n: v^Tw\le 0, \forall w \in \mathcal{N} \}$ is said to be the tangent cone to $\Omega$ at $u\in  \Omega.$
Further, the projection operator $P_\Omega(u)$ for $u\in \Omega$ is denoted as $P_{\varOmega}\left( u \right) =\arg \min _{v\in \varOmega}\lVert u-v \rVert $.
Moreover, for a direction $v$ and $u\in \Omega$, we define the differentiated projction operator as $
	\Pi_\Omega (u,v) =\lim_{\varepsilon \to 0} \frac{P_\Omega (u + \varepsilon v) -v}{\varepsilon}$ which is equivalent to $P_{T_\Omega(u)}(v)$. 
The calculation criteria of $\Pi_\Omega (u,v) $ is revealed in the following lemma, which is restated from \cite{weibao}.
\begin{lemma}\label{lem.qei3}
	The differentiated projection operator $\Pi_\Omega (u,v)$ can be computed according to the following properties:
(i)  if $u \in int(\Omega)$, then $\Pi_\Omega (u,v) =v;$  (ii)  $u \in \partial \Omega$ and $\max_{w \in n_\Omega (u)} \langle v, w \rangle \leq 0$, then $\Pi_\Omega (u,v) =v;$  
(iii)    $u \in \partial \Omega$ and $\max_{w \in n_\Omega (u)} \langle v, w \rangle \geq 0$, then $\Pi_\Omega (u,v) =v-\langle v, w^* \rangle w^*,$ where $w^*=\arg\max_{w \in n_\Omega (u)} \langle v, w \rangle$.
\end{lemma}

\subsection{Differential Inclusions and Nonsmooth Analysis}

By denoting a set-valued map  $\mathcal{F}: X  \rightrightarrows Y $, a differential inclusion is expressed as 
\begin{align}\label{di}
	\dot x \in \mathcal{F}(x), \quad x(0) = x_0.
\end{align}
The solution to (\ref{di}) in the caratheodory sence is an absolutely continuous function $x : [0, +\infty ) \rightarrow  \mathbb{R}^n$ which satisfies (\ref{di})  for almost all $t \in [0, +\infty ) $.
For a convex function $f : \Omega \rightarrow \mathbb{R}$ with $\Omega \subset \mathbb{R}^n$, $\partial f: \Omega \rightrightarrows \mathbb{R}^n$ is a nonempty, convex and compact set-valued map while being locally bounded and upper semicontinuous at $x\in \Omega$. Furthermore, if $\mathcal{F}$ in (\ref{di}) is a locally bounded and upper semicontinuous set-valued map with nonempty, convex and compact value, then there exists a solution to (\ref{di}).

To facilitate our analysis in subsequent, we review some results from \cite{Aubin1984Differential} as below.

\begin{lemma} \label{lem.set}
Consider a set-valued map $\mathcal{F}: X  \rightrightarrows Y $ and a closed convex subset $\Omega \subset \mathbb{R}^n$. For two differential inclusions 
	\begin{eqnarray}
	\label{eq:di1}
	&&\dot x \in \mathcal{F}(x) - \mathcal{N}_{\Omega}(x), \quad x(0) = x_0,\\
	\label{eq:di2}
	&&\dot x \in \Pi_{\Omega}(x,\mathcal{F}(x)), \quad x(0) = x_0,
\end{eqnarray}
we have the following satements hold.

1) There exists a solution to (\ref{eq:di1}) if $\mathcal{F}$ is bounded on $\Omega$.

2) The trajectory $x(t)$ is a solution to (\ref{eq:di1}) if and only if it is also a solution to (\ref{eq:di2}).
\end{lemma}

Denote a continuously differentiable function $V : \mathbb{R}^n  \rightarrow \mathbb{R}$. The set-valued Lie derivative for $V$ along with (\ref{di}) 
can be expressed as 
$$
\mathcal{L}_\mathcal{F}V = \{p \in \mathbb{R} \, | \, \exists \, v \in \mathcal{F}(x) {\text{ s.t. }}  v^T\nabla V(x) =p\}.
$$

The invariance principle with respect to differential inclusions is given below (see \cite{BG}).

\begin{lemma} \label{lem.LaSalle}
Suppose that $V: \mathbb{R}^n  \rightarrow \mathbb{R}$ is continuously differentiable and the compact subset $S \subset \mathbb{R}^n$ is a strongly  positively invariant set for (\ref{di}). The solution to (\ref{di}) starting from $S$ converges to the largest weakly positively invariant subset $M \subset S \cap \{x \in \mathbb{R}^n \, | \, 0 \in \mathcal{L}_\mathcal{F}V(x)\}$  when $\max \mathcal{L}_\mathcal{F}V \le 0$
or $\mathcal{L}_\mathcal{F}V=\emptyset $  for all $x \in S$, and the solution to (\ref{di}) is bounded.
\end{lemma}


\subsection{Problem Formulation}

The resource allocation problem considered in this paper is a constrained optimization problem over a network of $N$ agents. Each agent equips with a  local cost function (may be non-smooth) $f_i(x_i): \Omega_i \rightarrow  \mathbb{R}$, where $\Omega_i$ is a closed convex subset of $\mathbb{R}^n$, and $x_i \in \Omega_i$ is the local decision variable (or called an allocation). And all agents aim to cooperatively  minimize a global cost function $
f\left( x \right) =\sum_{i=1}^N{f_i\left( x_i \right)}$ through a communication digraph $\mathcal{G}$, where $x=col(x_1,\ldots,x_N) \in \mathbb{R}^{Nn}$.

Specifically, the nonsmooth resource allocation problem can be formulated as:
\begin{align} \label{op1}
	&\min_{x\in \mathbb{R}^{Nn}}f\left( x \right),   ~~ f\left( x \right) =\sum_{i=1}^N{f_i}\left( x_i \right) \notag \\ &\text {subject to}~ \sum_{i=1}^N x_{i} = \sum_{i=1}^N d_{i} \notag  \\&  x_{i} \in\Omega_{i}, ~~ i \in \{1,\dots,N\}
\end{align}
where $\sum_{j=1}^N d_{i}$ is the total network resource.

It follows from the constraints presented in (\ref{op1}) that all agents must satisfy resource allocation constraints $\sum_{i=1}^N x_{i} = \sum_{i=1}^N d_{i}$ and local feasibility contraints. Moreover, each agent $i$ only has the private knowledge about $f_i,d_i$ and $\Omega_{i}$.

The main task in this paper is to develop a distributed algorithm to seek the optimum for the nonsmooth resource allocation problem formulated in (\ref{op1}) under a weight-balanced digraph.

Some standard assumptions widely used in \cite{Deng2018,Zhu2019TSM,Deng2020}  are given below for the wellposedness of problem (\ref{op1}).

\begin{assumption}\label{ass.1}
	The communication digraph $\mathcal{G}$  is strongly connected and weight-balanced.
\end{assumption}

\begin{assumption}\label{ass.2}
	Slater's condition holds, that is, there exists an interior point $x_i \in int(\Omega_i)$, such that $\sum_{i=1}^N{x_i}=\sum_{i=1}^N{d_i}$, for $i \in \{1,\dots,N\}$.
\end{assumption}

\begin{assumption}\label{ass.3}
	The local cost function $f_i$ is $\omega$-strongly convex with constant $\omega > 0$, for $i \in \{1,\dots,N\}$.
\end{assumption}

\begin{remark}
	Problem (\ref{op1}) can be viewed as a generalized version of problems considered in \cite{Yi2016,Xuefang2019,Wangex2018,Chen2021TSM,2018Dis,YinAuto2019,Lgriad,He2017,Sunlanlan2021,KiaAuto2021,non108,non109,Lian22021,Zeng2016}. Specifically, compared to \cite{Yi2016,Xuefang2019,Wangex2018,Chen2021TSM,2018Dis,YinAuto2019,Lgriad,He2017,Sunlanlan2021,KiaAuto2021}, the assumption that the cost function is differentiable (or even twice differentiable) is removed. Compared to\cite{Yi2016,Xuefang2019,Chen2021TSM,non108,non109,Lian22021,Zeng2016}, undirected graph is extended to weight-balanced digraph. Compared to  \cite{2018Dis,Chen2021TSM,Lgriad,He2017}, the local box constraint is replaced by a more general form of convex set constraint.
\end{remark}

\begin{lemma} \label{lem.AA}
	Given a matrix $A\in \mathbb{R}^{n\times n}$, and a vector $x\in \mathbb{R}^n$, $A^TAx=0_n$  amounts to $Ax=0_n$.
\end{lemma}
Proof: It follows from $A^TAx=0_n$ that $(Ax)^T(Ax)=x^TA^TAx=0$. Then we can infer that $Ax=0_n$ by $(Ax)^T(Ax)=0$.\hfill $\Box$

To end this section, the optimality condition (see \cite{Deng2020}) associated with the problem (\ref{op1}) is given as below.
\begin{lemma} \label{lem.kkt}(see \cite[Theorem 3.34]{Franz2})
	For each $i \in \{1,\dots,N\}$, $\tilde{x}_i$ is an optimal solution of the problem (\ref{op1}), if there exists $\tilde{\mu}\in \mathbb{R}^n$ such that 
	\begin{equation}\label{kkt}
		\begin{split}
			& 0_n\in\partial f_i\left( \tilde{x}_i \right) -\tilde{\mu}+\mathcal{N}_{\varOmega _i}\left( \tilde{x}_i \right)\\
			&\sum_{i=1}^N{\tilde{x}_i}=\sum_{i=1}^N{d_i}.
		\end{split}
	\end{equation}
Conversely, for each $i \in \{1,\ldots,N\}$, there exists $\tilde{\mu} \in \mathbb{R}^n$ such that the condition (\ref{kkt}) is satisfied, if  $\tilde{x}_i$ is an optimal solution of problem (\ref{op1}).
\end{lemma}


\section{Main Results}\label{sec.r}
In this section, to tackle the nonsmooth resource allocation problem, we first design an distributed algorithm with considering local feasibility constraints, and then we design a simplified algorithm for solving the differentiable version of the problem without local constraints in Section \ref{A}. Afterwards, we analyze their convergence in Section \ref{B}.
\subsection{Distributed Algorithm Design}\label{A}
In this subsection, we first provide a distributed algorithm equipped with the differentiated projection operator for the problem (\ref{op1}). Then, for the possibility of exponential convergence, we design the other algorithm for a simplified problem (\ref{op2}) by imposing the differentiability of cost function and discarding the local feasibility constraints.

Motivated by \cite{Yi2016,Deng2018}, differentiated projection operations coupled with differential inclusions are utilized to design the following algorithm for agent $i \in \{1,\ldots,N\}$, aiming to tackle the problem (\ref{op1}):
\begin{align} \label{al1}
	\begin{cases}             
		\dot{x}_i\in \varPi _{\varOmega _i}\left( x_i,-\partial f_i\left( x_i \right) +\mu _i \right) \quad x_i(0) \in \Omega_i \\
		\dot{\mu}_i=k_1\left( \eta _i-x_i+d_i \right) -k_2\sum_{j=1}^N{a_{ij}\left( \mu _i-\mu _j \right)}\\
		\dot{\eta}_i=-k_3\sum_{j=1}^N{a_{ij}\left( \left( \eta _i-x_i+d_i \right)  -\left( \eta _j-x_j+d_j \right)\right)}
	\end{cases}
\end{align}


where \begin{align} \label{can1}
	k_1>\frac{\lVert L \rVert ^2}{\hat{\lambda}_2\omega}, \quad k_2>\frac{k_{1}^{2}}{\hat{\lambda}_{2}^{2}},  \quad  k_3>0.
\end{align}

Furthemore, by dropping the local contrains, the problem (\ref{op1}) can be simply described as below:
\begin{align} \label{op2}
	\min_{x\in \mathbb{R}^{Nn}}&~~f\left( x \right) ,\quad f\left( x \right) =\sum_{i=1}^N{f_i}\left( x_i \right) \notag \\ &\text {subject to}~ \sum_{i=1}^N x_{i} = \sum_{i=1}^N d_{i}.
\end{align}

When encountering the case that cost functions are differentiable with $\theta$-Lipschitz gradients, the algorithm (\ref{al1}) reduce to 
\begin{align} \label{al2}
	\begin{cases}           
		\dot{x}_i=-\nabla f_i\left( x_i \right) +\mu _i\\
		\dot{\mu}_i=k_1\left( \eta _i-x_i+d_i \right) -k_2\sum_{j=1}^N{a_{ij}\left( \mu _i-\mu _j \right)}\\
		\dot{\eta}_i=-k_3\sum_{j=1}^N{a_{ij}\left( \left( \eta _i-x_i+d_i \right) -\left( \eta _j-x_j+d_j \right) \right)}
	\end{cases}
\end{align}
where 
\begin{align} \label{can2}
	k_1>\max \left\{ \frac{\lVert L \rVert ^2\left( \omega +1 \right)}{\hat{\lambda}_2\omega ^2}+\frac{\theta ^2}{2w},1 \right\}, \quad k_2>\frac{k_{1}^{2}}{\hat{\lambda}_{2}^{2}}, \quad k_3>0.
\end{align}

\begin{remark}
The above algorithms (\ref{al1}) and (\ref{al2}) have the same structure based on the gradient descent scheme to seek the optimum. Specifically, each agent has three state variables $x_i,\mu_{i}$, and $\eta_i$, where $x_i$ is the decision variable that allows each agent to seek the optimal solution that minimizes the overall cost function with the help of the auxiliary variables $\mu_{i}$ and $\eta_i$. Moreover, together with $\mu_{i}$, $\eta_i$ can be utilized to satisfy the coupling equality constraint in problem (\ref{op1}).
Note that the local set constraint is handled in algorithm (\ref{al1}) using the differential projection operator, which can be obtained by differentiating the projection operator in \cite{Zhu2019TSM,Deng2020,Zeng2016}, and in this way, the decision variables are guaranteed to always lie within the constraint set.
\end{remark}
\begin{remark}
It is necessary to clarify that although private information $x_j$ and $d_j$ of neighbors appear in the dynamics of agents, the information that actually interacts among agents is the overall information $\eta_j+x_j-d_j$. Due to the private properties of $\eta_j$, $x_j$ and $d_j$, each agent cannot distinguish the specific values of the others' decision information $x_j$ or local resource $d_j$ based on the sum information $\eta_j+x_j-d_j$.  Therefore, algorithms (\ref{al1}) and (\ref{al2}) still protect privacy to some extent.
and do not add additional communication burden to the system than using the algorithm  in \cite{Yi2016,Deng2018,Zhu2019TSM}.
\end{remark}
\begin{remark}
It will be shown in the subsequent proofs that the parameter ranges in algorithm (\ref{al1}) and algorithm (\ref{al2}) are necessary to guarantee the convergence of  algorithms. Although the determination of the parameter ranges requires the values of $\theta, \omega, \hat{\lambda}_2$  and  $\lVert L \rVert$  in a similar way as in \cite{Deng2018,Deng2020}, we can use the distributed algorithms developed in existing literature (see \cite{concensus}) to calculate or estimate them up front so as to determine the parameter ranges, for enabling a fully distributed implementation of algorithms.
\end{remark}

\subsection{Convergence Analysis}  \label{B}
In this subsection, we investigate the correctness and convergence of  the two algorithms, respectively.

To begin with, we study the equilibrium point of (\ref{al1}), and then investigate the convergence of  (\ref{al1}).

Let
\begin{align*}
	x&=col\left( x_1,\cdots ,x_N \right)\\
	\mu &=col\left( \mu _1,\cdots ,\mu _N \right) 
	\\
	\eta &=col\left( \eta _1,\cdots ,\eta _N \right)\\
	d&=col\left( d_1,\cdots ,d_N \right) \\
	\varOmega &=\varOmega _1\times \cdots \times \varOmega _N .
\end{align*}

For convenience, we use the equivalent compact form replacing (\ref{al1}) as  follows:
\begin{align} \label{al11}
	\begin{cases}            
		\dot{x}\in \varPi _{\varOmega}\left( x,-\partial f\left( x \right) +\mu \right)  \quad x(0) \in \Omega \\
		\dot{\mu}=k_1\left( \eta -x+d \right) -k_2\left( L\otimes I_n \right) \mu\\ 
		\dot{\eta}=-k_3\left( L\otimes I_n \right) \left( \eta -x+d \right) .
	\end{cases}
\end{align}
Obviously, Lemma  \ref{lem.set} can guarantee that there exists a solution to (\ref{al11}).

The next lemma given below reveals the relationship between the equilibrium  of (\ref{al1}) and the optimum of (\ref{op1}).
	\begin{lemma}\label{lem1}
Suppose that Assumption \ref{ass.1} holds. If the initial condition $\sum_{i=1}^N{\eta_i(0)}=0_n$ is satisfied and $(x^*,\mu^*,\eta^*)$ is an equilibrium point of algorithm (\ref{al11}), then $x^*$ is an optimal solution of  problem (\ref{op1}). Conversely, if $x^*$ is an optimal solution of  problem (\ref{op1}), then there exists $( \mu^*,\eta^*)\in \mathbb{R}^{Nn}\times \mathbb{R}^{Nn}$ such that $(x^*,  \mu^*, \eta^*)$ is an equilibrium point of algorithm (\ref{al11}).
\end{lemma}
\emph{Proof:} 
1) For a given equilibrium point $\left( x^*,\mu ^*,\eta ^* \right)$ of (\ref{al11}), it is verified that
\begin{subequations}
	\begin{align}    \label{Th1a}        
		0_{Nn}\in& \varPi _{\varOmega}\left( x^*,-\partial f\left( x^* \right) +\mu ^* \right) \\ \label{Th1b}
		0_{Nn}=&k_1\left( \eta ^*-x^*+d \right) -k_2\left( L\otimes I_n \right) \mu ^*\\  \label{Th1c}
		0_{Nn}=&-k_3\left( L\otimes I_n \right) \left( \eta ^*-x^*+d \right).
	\end{align}
\end{subequations}

In the light of Lemma \ref{lem.qei3}, we have $\dot{x}_i\in \ T_{\varOmega _i}\left( x_i \right) $.  Knowing that $
\dot{x}_i\left( 0 \right) \in \varOmega _i $, using the viability theorem of \cite{Aubin1984Differential}, $
{x}_i\left( t\right) \in \varOmega _i $ holds for all $t \ge 0$.

By virtue of the strong connectedness  and weight-balance of digraphs, there exists $\vartheta \in \mathbb{R}^n$ suth that 
\begin{align} \label{tha} 
	(L\otimes I_n)\mu^*=1_N\otimes \vartheta
\end{align}
based on (\ref{Th1b}) and (\ref{Th1c}).

Multiplying left by $L^T\otimes I_n$ on bath sides of (\ref{tha}), it obtains $(L^TL\otimes I_n)\mu^*=L^T1_N \otimes \vartheta=0_{Nn}$.
Further, invoking Lemma \ref{lem.AA}, it implies $(L\otimes I_n)\mu^*=0_{Nn}$
which amounts  to $\mu_i^*=\mu_j^*, \  \forall i, j \in \{1,\dots,N\}$.

On the other hand, it follows from (\ref{Th1b}) that $k_1\left( 1_{N}^{T}\otimes I_n \right) \left( \eta ^*-x^*+d \right) =k_2\left( 1_{N}^{T}L\otimes I_n \right) \mu ^*=0_{n}$. Meanwhile, it can be seen that $(1_N^T\otimes I_n)\dot{\eta}(t)=0_{n}$ by (\ref{al11}). As a result, with the given initial value $(1_N^T\otimes I_n)\eta(0)=0_{n}$,  we know that $(1_N^T\otimes I_n)\eta(t)=0_{n}.$ 
Hence, we can deduce that $k_1(1_N^T\otimes I_n)(-x^*+d)=0_n$, i.e. $\sum_{i=1}^{N}x_i^*=\sum_{i=1}^{N}d_i$.

Moreover, on the basis of Lemma \ref{lem.qei3}, one of the following cases occurs derived from (\ref{Th1a}).
\begin{itemize}
\item[(i)] $x_i^* \in \partial \Omega_{i}$ and $0_n \in -\partial f_i(x_i^*)+\mu_{i}^*$.
\item[(ii)] $x_i^* \in \partial \Omega_{i}$ and $ -\partial f_i(x_i^*)+\mu_{i}^* \in \mathcal{N}_{\Omega_{i}}(x_i^*)$.
\item[(iii)] $x_i^* \in int(\Omega_{i}) $  and $0_n \in -\partial f_i(x_i^*)+\mu_{i}^*+\mathcal{N}_{\Omega_{i}}(x_i^*)$.
\end{itemize}
Obviously, it maniffests that $0_n \in \partial f_i(x_i^*)-\mu_i^*+\mathcal{N}_{\Omega_i}(x_i^*)$.

Hence, $x^*$ is an optimal solution of  (\ref{op1}) based on Lemma \ref{lem.kkt}.

2) For an optimal solution $x^*$ of  (\ref{op1}), it leads to the following conditions hold:
\begin{subequations}\label{zui2} 
	\begin{align}  \label{zui2a}         
		&0_{n}\in \partial f_i\left( x_{i}^{\ast} \right) -\mu _{i}^{\ast}+\mathcal{N}_{\varOmega _i}\left( x_{i}^{\ast} \right) \quad  \\ \label{zui2b}     
&\sum_{i=1}^N{x_{i}^{\ast}}=\sum_{i=1}^N{d_i}\\\label{zui2c}
&\mu _{i}^{\ast}=\mu _{j}^{\ast},\   \forall i,j\in \{1,...,N\}\\ \label{zui2d}
&x_{i}^{\ast}\in \varOmega _i, \  \forall i\in \{1,...,N\}.
	\end{align}
\end{subequations}

We can obtain $(L\otimes I_n)\mu^*=0_{Nn}$ by taking $\mu^*=col\{\mu_1^*,\dots,\mu_N^* \}$, on account of (\ref{zui2c}).

On the grounds of (\ref{zui2a}),  it can be seen that (\ref{Th1a}) is satisfied. Further, consider $\eta^*=col\{\eta_1^*,\dots,\eta_N^*\}$ defined by $\eta_i^*=x_i^*-d_i$. It gives rise to that (\ref{Th1b}) and (\ref{Th1c}) can be established. \hfill $\Box$

Subsequently, having the above discussion, we can see that $(x^*,\mu^*,\eta^*)$ is an equilibrium point of (\ref{al11}).

By Lemma \ref{lem1}, we have the following statement corresponding to the convergence of (\ref{al11}).
	\begin{theorem}\label{the1}
	Suppose that Assumptions \ref{ass.1}, \ref{ass.2}, and \ref{ass.3} hold. The trajectory $(x(t),\mu (t),\eta (t))$ starting from any initial state $(x(0), \mu(0), \eta(0))$ satisfying $ \sum_{i=1}^N{\eta _i\left( 0 \right)}=0_n $ under (\ref{al11}) converges asymptotically to the equilibrium of algorithm (\ref{al11}), and $x(t)$ converges asymptotically to the optimum of  problem (\ref{op1}).
\end{theorem}
\emph{Proof:} 
For notational simplicity, we assume $n=1$ without loss of generality.
Before moving on, in view of Lemma \ref{lem.qei3}, it is clear that 
$$
\varPi _{\varOmega _i}\left( x_i,-\partial f_i\left( x_i \right) +\mu _i \right) \subset -\partial f_i\left( x_i \right) +\mu _i-\mathcal{N}_{\varOmega _i}\left( x_i \right) 
$$
and thus $-\partial f_i\left( x_{i}^{*} \right) +\mu _{i}^{*}\subset \mathcal{N}_{\varOmega _i}\left( x_{i}^{*} \right) $.

To clarify the subsequent analysis, we convert (\ref{al11}) into a concise form. In detail, we first transform the equilibrium point of system (\ref{al11})  to the origin through the coordinate transformation as below:
\begin{align*}
\bar{x}&=x-x^*\\
\bar{\mu}&=\mu -\mu ^*\\
\bar{\eta}&=\eta -\eta ^*.
\end{align*}

Then we have the following system which is identical to (\ref{al11}):
\begin{equation} 
	\begin{cases} \dot {\bar x} \in \bigg \{ p \in \mathbb {R}^{N}|~p= - \varrho + \bar \mu \\ \qquad \varrho \in \partial f(x) - \partial f({x^{*}}) +\mathcal{N} _{\Omega}(x) - \mathcal{N}_{\Omega}\left ({x^{*}}\right )\bigg \} \\ 
		\dot{\bar{\mu}}=k_1\left( \bar{\eta}-\bar{x} \right) -k_2L\bar{\mu}\\ 
		\dot{\bar{\eta}}=-k_3L\left( \bar{\eta}-\bar{x} \right) . \end{cases} 
\end{equation}

Further, by using the following orthogonal transformation introduced in Lemma \ref{lem.tu} as below:
\begin{equation}\label{zhenjiao}
	\begin{split}
		\varepsilon &=col\left( \varepsilon _1,\varepsilon _2 \right) =\left[ \begin{matrix}
			r&		R\\
		\end{matrix} \right] ^T\bar{x}\\
		\xi &=col\left( \xi _1,\xi _2 \right) =\left[ \begin{matrix}
			r&		R\\
		\end{matrix} \right] ^T\bar{\mu}\\
		\zeta &=col\left( \zeta _1,\zeta _2 \right) =\left[ \begin{matrix}
			r&		R\\
		\end{matrix} \right] ^T\bar{\eta}
	\end{split}
\end{equation}
where $\varepsilon _1,\xi _1,\zeta _1 \in \mathbb{R}  $ and  $\varepsilon _2,\xi _2,\zeta _2 \in \mathbb{R}^{N-1}  $, we recast (\ref{al11}) as
\begin{subequations} \label{zhen}
	\begin{eqnarray}
	&\begin{cases} 
		{\dot \varepsilon }_{1} \in \bigg \{p \in \mathbb {R}|~p={\xi _{1}} - {r^{\mathrm{ T}}}\varrho \\ \quad \qquad  \varrho \in \partial f(x) - \partial f({x^{*}}) + \mathcal{N}_{\Omega}(x) - \mathcal{N}_{\Omega}\left ({x^{*}}\right )\bigg \} \\ 
	\dot{\xi}_1=-k_1\varepsilon _1 \\ 
		\dot{\zeta}_1=0 \end{cases}  
	\\
	&\begin{cases} {\dot \varepsilon _{2}} \in \bigg \{p \in \mathbb {R}^{N-1}|~p=  {{\xi _{2}} - {R^{\mathrm{ T}}} \varrho}  \\ \quad \qquad  \varrho \in \partial f(x) - \partial f\left ({{x^{*}}}\right ) + \mathcal{N}_{\Omega}(x) - \mathcal{N}_{\Omega}\left ({x^{*}}\right )\bigg \} \\ 
	\dot{\xi}_2=-k_1\varepsilon _2+k_1\zeta _2-k_2R^TLR\xi _2\\ 
	\dot{\zeta}_2=-k_3R^TLR\left( \zeta _2-\varepsilon _2 \right).   \end{cases} 
	\end{eqnarray}
\end{subequations}

According to the above arguments, it is obvious that $x$ is convergent to the optimal solution of problem (\ref{op1}) if $\varepsilon$ converges toward the origin. In what follows, we study the convergence of (\ref{zhen}).

Construct the Lyapunov function candidate as 
\begin{align}   
	V_1=&\frac{k_1}{2}\left( \lVert \varepsilon _1 \rVert ^2+\lVert \varepsilon _2 \rVert ^2 \right) +\frac{1}{2}\left( \lVert \xi _1 \rVert ^2+\lVert \xi _2 \rVert ^2 \right)  \nonumber \\
	&+\frac{1}{2k_3}\lVert \zeta _2 \rVert ^2
\end{align} 
where $k_1$ and $k_3$ satisfies (\ref{can1}).

The simple calculations can give rise to the set-valued Lie derivative of $V_1$ with respect to (\ref{zhen}) as below:
\begin{align} 
	{\mathcal {L}_{(\ref{zhen})} V}_{1}=&\bigg \{ p \in \mathbb {R}|~p=-k_1\left( \varepsilon _{1}^{T}r^T\varrho+\varepsilon _{2}^{T}R^T\varrho \right) -k_1\xi _{2}^{T}\zeta _2\nonumber \\
	&\quad -k_2\xi _{2}^{T}R^TSym\left( L \right) R\xi _2-\zeta _{2}^{T}R^TLR\varepsilon _2  \nonumber \\
	&\quad-\zeta _{2}^{T}R^TSym\left( L \right) R\zeta _2  \nonumber \\
	&\quad \varrho \in \partial f(x) - \partial f({x^{*}}) + \mathcal{N}_{\Omega}(x) - \mathcal{N}_{\Omega}\left ({x^{*}}\right )\bigg \}.\qquad 
\end{align}

Note that 
\begin{align*} \left \langle{ {x_{i} - x_{i}^{*},{\mathcal{N}_{\Omega _{i}}}(x_{i})} }\right \rangle\subseteq&\mathbb {R}_{+} \\ \left \langle{ {x_{i} - x_{i}^{*},{\mathcal{N}_{\Omega _{i}}}\left ({x_{i}^{*}}\right )} }\right \rangle\subseteq&\mathbb {R}_{-} ,
\end{align*}
it signifies 
\begin{equation*}
	 \langle \bar x,{N_{\Omega}}(x) - {\mathcal{N}_{\Omega}}(x^{*})\rangle \subseteq \mathbb {R}_{+}. 
\end{equation*}
Besides, it is clear that 
\begin{align*} 
	-\varepsilon _{1}^{\mathrm{ T}}{r^{\mathrm{ T}}}\varrho - \varepsilon _{2}^{\mathrm{ T}}{R^{\mathrm{ T}}} \varrho \in&- \left \langle{ \bar x, \partial f(x) - \partial f\left ({{x^{*}}}\right )}\right \rangle \\&-\left \langle{ \bar x,{\mathcal{N}_{\Omega}}(x_{i}) - {\mathcal{N}_{\Omega}}\left ({x_{i}^{*}}\right )}\right \rangle . 
\end{align*}
Thus, in light of the orthogonal transformation (\ref{zhenjiao}) and the strong convexity of cost function, we can establish the inequality as follows:
\begin{align} \label{20}
	&\max \bigg \{p \in \mathbb {R} |~p= - k_{1} \left ({\varepsilon_{1}^{\mathrm{ T}}{r^{\mathrm{ T}}}\varrho + \varepsilon _{2}^{\mathrm{ T}}{R^{\mathrm{ T}}}\varrho}\right ) \nonumber\\ 
	&\qquad\;\;~ \varrho \in \partial f(x) - \partial f\left ({{x^{*}}}\right ) + \mathcal{N}_{\Omega}(x) - \mathcal{N}_{\Omega}\left ({x^{*}}\right )\bigg \}  \nonumber \\
	& \qquad  \;\;  \leq -k_{1} \omega {\left \|{\varepsilon }\right \|^{2}}. 
\end{align}

 By using the Yong's inequality, simple computations lead to 
\begin{align}   
	-k_1\xi _{2}^{T}\zeta _2\le \frac{k_{1}^{2}}{\hat{\lambda}_2}\lVert \xi _2 \rVert ^2+\frac{\hat{\lambda}_2}{4}\lVert \zeta _2 \rVert ^2,
\end{align} 
\begin{align}   
	-\zeta _{2}^{T}R^TLR\varepsilon _2\le \frac{\hat{\lambda}_2}{4}\lVert \zeta _2 \rVert ^2+\frac{\lVert L \rVert ^2}{\hat{\lambda}_2}\lVert \varepsilon \rVert ^2.
\end{align} 

Moreover, by considering the fact that the communication graph is weight-balanced and strongly connected, one has that
\begin{align}  \label{L1}
	-k_2\xi _{2}^{T}R^TSym\left( L \right) R\xi _2\le -k_2\hat{\lambda}_2\lVert \xi _2 \rVert ^2,
\end{align} 
\begin{align}  \label{L2} 
	-\zeta _{2}^{T}R^TSym\left( L \right) R\zeta _2\le -\hat{\lambda}_2\lVert \zeta _2 \rVert ^2.
\end{align} 

By resorting to the above inequalities (\ref{20})-(\ref{L2}), simple manipulations yeild
\begin{align}  \label{xita1}
\max \mathcal{L}_{(\ref{zhen})}\le -\tau
\end{align} 
where 
\begin{align}   \label{xita2}
	\tau=&\left( k_1\omega -\frac{\lVert L \rVert ^2}{\hat{\lambda}_2} \right) \lVert \varepsilon \rVert ^2\nonumber \\
	&+\left( k_2\hat{\lambda}_2-\frac{k_{1}^{2}}{\hat{\lambda}_2} \right) \lVert \xi _2 \rVert ^2\nonumber \\
	&+\frac{\hat{\lambda}_2}{2}\lVert \zeta _2 \rVert ^2 \ge 0.
\end{align} 
Equipped with (\ref{xita1}) and (\ref{xita2}), we can deduce that $\varepsilon_1,\varepsilon_2,\xi_1,\xi_2$ and $\zeta_2$ are bounded.

Invoking Lemma \ref{lem.LaSalle}, it  verifies that the solution of (\ref{zhen}) is convergent to the following set $M$:
\begin{align*} 
	M=&\big \{\varepsilon \in {\mathbb {R}^{N}}, \xi _{2} \in {\mathbb {R}^{N-1}} \\&~~ {\zeta _{2}} \in {\mathbb {R}^{N - 1}} | \varepsilon = 0_{N},\xi _{2} = 0_{N-1},{\zeta _{2}} = 0_{N-1} \big \}. 
\end{align*}
Furthermore, on the grounds of Lemma \ref{lem1}, we can obtain that  $\lim _{t\rightarrow \infty}x\left( t \right) =x^*$    
without much effort, which can complete the proof.\hfill $\Box$

In what follows, we depict the property with respect to the equilibrium of (\ref{al2}), and then analyze the exponential convergence of (\ref{al2}). We rewrite (\ref{al2}) in a compact form as :
\begin{align} \label{al22}
	\begin{cases}             
		\dot{x}=-\nabla f\left( x \right) +\mu \\
		\dot{\mu}=k_1\left( \eta -x+d \right) -k_2\left( L\otimes I_n \right) \mu \\
		\dot{\eta}=-k_3\left( L\otimes I_n \right) \left( \eta -x+d \right) .
	\end{cases}
\end{align}

Similar to Lemma \ref{lem1}, the following result is concerning the equilibrium  of (\ref{al22}).
	\begin{lemma}\label{lem2}
	Suppose that Assumption \ref{ass.1}  holds. With the initial state satisfying $\sum_{i=1}^N{\eta_i(0)}=0_n$, if $(x^*,\mu^*,\eta^*)$ is an equilibrium of algorithm (\ref{al22}), then $x^*$ is an optimal solution of  problem (\ref{op2}). Converely, if $x^*$ is an optimum of problem (\ref{op2}), then there exists $( \mu^*,\eta^*)\in \mathbb{R}^{Nn}\times \mathbb{R}^{Nn}$ such that $(x^*,  \mu^*,\eta^*)$ is an equilibrium of algorithm (\ref{al22}).
\end{lemma}
\emph{Proof:} 
1) Consider an equilbrium $(x^*,\mu^*,\eta^*)$ of (\ref{al22}), one can obtain
\begin{subequations} 
	\begin{align}   \label{Th2a}
		0_{Nn}=&-\nabla f\left( x^* \right) +\mu ^* \\ \label{Th2b}
		0_{Nn}=&k_1\left( \eta ^*-x^*+d \right) -k_2\left( L\otimes I_n \right) \mu ^* \\ \label{Th2c}
		0_{Nn}=&-k_3\left( L\otimes I_n \right) \left( \eta ^*-x^*+d \right) .
	\end{align} 
\end{subequations}

%

Following the similar demonstration as in Lemma \ref{lem.AA}, we can see that $\sum_{i=1}^N{x_{i}^{*}}=\sum_{i=1}^N{d_i}$ and 
$\mu _{i}^{*}=\mu _{j}^{*},\ \forall \ i,j \in \{1,\dots,N\}$ are valid.

Moreover, it is obvious that $\nabla f_i\left( x_i^* \right)=\nabla f_j\left( x_j^* \right),  \forall \ i,j \in \{1,\dots,N\}$ from (\ref{Th2a}) and $
\nabla f\left( x^* \right) =col\left( \nabla f_1\left( x_{1}^{*} \right) ,\ldots ,\nabla f_N\left( x_{N}^{*} \right) \right) 
$.

Via the previous arguments, according to \cite[Lemma 3.1]{siam}, it is obtained that $x^*$ is an optimum of  (\ref{op2}).

2) If $x^*$ is an optimal solution of  (\ref{op2}), then there exists $\mu_i^* \in \mathbb{R}^n$ such that 
\begin{align} \label{kt}
	\begin{cases}             
\nabla f_i\left( x_{i}^{*} \right) =\nabla f_j\left( x_{j}^{*} \right),  \ \forall \ i,j\in \left\{ 1,\cdots ,N \right\} \\
\sum_{i=1}^N{x_{i}^{*}}=\sum_{i=1}^N{d_i}\\
\mu _{i}^{*}=\mu _{j}^{*}, \ \ \forall \ i,j\in \left\{ 1,\cdots ,N \right\}.
	\end{cases}
\end{align}

Let $\mu^*=col\{\mu_1^*,\dots,\mu_N^*\}$.  Then (\ref{kt}) indicates $(L\otimes I_n)\mu^*=0_{Nn}$.
Furthermore, by choosing $\eta^*=x^*-d$, it can be claimed that (\ref{Th2b}) and (\ref{Th2c}) hold.

Thus, $(x^*,\mu^*,\eta^*)$ is an equilibrium of (\ref{al22}).\hfill $\Box$

Next, we discuss the exponential convergence property of (\ref{al22}).
\begin{theorem}\label{the2}
Suppose that Assumptions \ref{ass.1},  \ref{ass.2} and \ref{ass.3} hold.
If the gradients of local cost functions are $\theta$-Lipschitz continuous,
then  algorithm (\ref{al22}) with initial state satisfying $ \sum_{i=1}^N{\eta _i\left( 0 \right)}=0_n $ can exponentially converge to the optimum of  problem (\ref{op2}).
\end{theorem}
\emph{Proof:} 
Without loss of generality, we assume $n=1$ throughout the following analysis. Employing a similar manipulations in Theorem \ref{the1}, we study the convergence of the following system as:
\begin{subequations} \label{zhen2}
	\begin{eqnarray}
	&\begin{cases} 
	{\dot \varepsilon }_{1} = {\xi _{1}} - {r^{\mathrm{ T}}}\varrho  \\ 
	\dot{\xi}_1=-k_1\varepsilon _1 \\ 
	\dot{\zeta}_1=0  
	\end{cases} \quad \quad \quad \quad \quad \quad  \quad \quad\\
	&\begin{cases} 
	      \dot{\varepsilon}_2=	{{\xi _{2}} - {R^{\mathrm{ T}}} \varrho}\\
				\dot{\xi}_2=-k_1\varepsilon _2+k_1\zeta _2-k_2R^TLR\xi _2\\ 
	\dot{\zeta}_2=-k_3R^TLR\left( \zeta _2-\varepsilon _2 \right) 
	\end{cases} 
	\end{eqnarray}
\end{subequations}
where $\varrho=\nabla f\left( x \right) -\nabla f\left( x^* \right) $.

Note that if $\varepsilon$ converges to 0, we can conclude that $x$ is convergent to the optimum of problem (\ref{op2}).

With $k_1$ given in (\ref{can2}), choose the following Lyapunov function candidate as
\begin{align}   
	V_2=&\frac{1}{2}\left( \frac{\omega +1}{\omega}k_1-1 \right) \lVert \varepsilon _1 \rVert ^2+\frac{1}{2\omega}\lVert \xi _1 \rVert ^2  \nonumber\\
	&+\frac{\omega +1}{2\omega}\left( \lVert \varepsilon _2 \rVert ^2+\lVert \xi _2 \rVert ^2+\lVert \zeta _2 \rVert ^2 \right) \nonumber\\
	&+\frac{1}{2}\lVert \varepsilon _1-\xi _1 \rVert ^2 .
\end{align} 
The derivative of $V_2$ along (\ref{zhen2}) is deduced drectly as 
\begin{align}   
	\dot{V}_2=&-\frac{\omega +1}{\omega}k_1\varepsilon _{1}^{T}r^T\varrho-\frac{\omega +1}{\omega}k_1\varepsilon _{2}^{T}R^T\varrho-\lVert \xi _1 \rVert ^2  \nonumber\\
	&+k_1\lVert \varepsilon _1 \rVert ^2-\frac{\omega +1}{\omega}k_2\xi _{2}^{T}R^TSym\left( L \right) R\xi _2\nonumber\\
	&-\frac{\omega +1}{\omega}k_1\xi _{2}^{T}\zeta _2 
	-\frac{\omega +1}{\omega}\zeta _{2}^{T}R^TSym\left( L \right) R\zeta _2  \nonumber\\
	&-\xi _{1}^{T}r^T\varrho+\frac{\omega +1}{\omega}\zeta _{2}^{T}R^TLR\varepsilon _2.
\end{align}

The Young's inequality and the Lipschitz property of the gradients of cost functions can be leveraged for getting
\begin{align}  \label{lies}
	-\xi _{1}^{T}r^T\varrho\le \frac{\theta ^2}{2}\lVert \varepsilon \rVert ^2+\frac{1}{2}\lVert \xi_1 \rVert ^2.
\end{align}

The orthogonal transformation (\ref{zhenjiao}) and the strong convexity of cost functions can be used for inferring the inequdity as follows
\begin{align}   
-\frac{\omega +1}{\omega}k_1\left( \varepsilon _{1}^{T}r^T\varrho+\varepsilon _{2}^{T}R^T\varrho \right) \le -k_1\left( \omega +1 \right) \lVert \varepsilon \rVert ^2.
\end{align}

In the same way as (\ref{L1}) and (\ref{L2}), we can obtain
\begin{align}   
	\frac{\omega +1}{\omega}k_2\xi _{2}^{T}R^TSym\left( L \right) R\xi _2\ge \frac{k_2\hat{\lambda}_2\left( \omega +1 \right)}{\omega}\lVert \xi _2 \rVert ^2,
\end{align} 
\begin{align}   
	\frac{\omega +1}{\omega}\zeta _{2}^{T}R^TSym\left( L \right) R\zeta _2\ge \frac{\hat{\lambda}_2\left( \omega +1 \right)}{\omega}\lVert \zeta _2 \rVert ^2.
\end{align}

Additionally, we can make use of the Young's inequality for obtaining
\begin{align}   
	\frac{\omega +1}{\omega}\zeta _{2}^{T}R^TLR\varepsilon _2\le& \frac{\hat{\lambda}_2\left( \omega +1 \right)}{4\omega}\lVert \zeta _2 \rVert ^2  \nonumber\\
	&+\frac{\lVert L \rVert ^2\left( \omega +1 \right)}{\hat{\lambda}_2\omega}\lVert \varepsilon _2 \rVert ^2,
\end{align} 
\begin{align}   \label{nice}
	-\frac{\omega +1}{\omega}k_2\xi _{2}^{T}\zeta _2\le \frac{k_{2}^{2}\left( \omega +1 \right)}{\hat{\lambda}_2\omega}\lVert \xi _2 \rVert ^2+\frac{\hat{\lambda}_2\left( \omega +1 \right)}{4\omega}\lVert \zeta _2 \rVert ^2.
\end{align} 

Subsequently, with (\ref{lies})-(\ref{nice}) in place, it implies

\begin{align}   \label{v2dao}
	\dot{V}_2\le& -\left( k_1\omega -\frac{\lVert L \rVert ^2\left( \omega +1 \right)}{\hat{\lambda}_2\omega}-\frac{\theta ^2}{2} \right) \lVert \varepsilon \rVert ^2-\frac{1}{2}\lVert \xi _1 \rVert ^2 \nonumber\\
	&-\frac{\omega +1}{\omega}\left( k_2\hat{\lambda}_2-\frac{k_{1}^{2}}{\hat{\lambda}_2} \right) \lVert \xi _2 \rVert ^2 \nonumber\\
	&-\frac{\hat{\lambda}_2\left( \omega +1 \right)}{2\omega}\lVert \zeta _2 \rVert ^2 \le -\rho \lVert \psi \lVert^2
\end{align} 
where 
$\psi =col\left( \varepsilon, \xi _1,\xi _2,\zeta _2 \right) $,  
$\rho =\min \left\{ \rho_1,\rho_2,\rho_3,\frac{1}{2}\right\} $ , $\rho_1=k_1\omega -\frac{\lVert L \rVert ^2\left( \omega +1 \right)}{\hat{\lambda}_2\omega}-\frac{\theta ^2}{2}$, $\rho_2=\frac{\omega +1}{\omega}\left( k_2\hat{\lambda}_2-\frac{k_{1}^{2}}{\hat{\lambda}_2} \right) $, and $\rho_3=\frac{\hat{\lambda}_2\left( \omega +1 \right)}{2\omega}$.

Notably, $V_2$ can be recast into $V_2=\frac{1}{2}\psi^T\varPhi \psi$, where $\varPhi=\text{diag}\{\varPhi_1,\varPhi_2\}$ with
\begin{align*} 
\Phi _1=& \begin{bmatrix}
	k\frac{\omega +1}{\omega}&		0_{N-1}^{T}&		-1\\
	0_{N-1}&		\frac{\omega +1}{\omega}I_{N-1}&		0_{N-1}\\
	-1&		0_{N-1}^{T}&		\frac{\omega +1}{\omega}\\
\end{bmatrix} \\
\Phi _2=&\begin{bmatrix}
	\frac{\omega +1}{\omega}I_{N-1}&		0_{\left( N-1 \right) \times \left( N-1 \right)}\\
0_{\left( N-1 \right) \times \left( N-1 \right)}&		\frac{\omega +1}{\omega}I_{N-1}\\
\end{bmatrix}. \quad  \quad \quad  \quad   \quad  \quad 
\end{align*}

To proceed, with reference to (\ref{v2dao}) and using the maximum and minimum eigenvalues of $ \varPhi$ defined by $\mu_{\max}$ and $\mu_{\min}$, respectively, it gives rise to
$$\frac{1}{2}\mu _{\min}\lVert \psi \rVert ^2\le V_2 \le \frac{1}{2}\mu _{\max}\lVert \psi \rVert ^2$$
which signifies 
\begin{align} \label{biiao}
	-\rho \lVert \psi \rVert ^2 \le -\frac{2\rho}{\mu_{\max}}V_2.
\end{align}

On account of (\ref{v2dao}) and (\ref{biiao}), simple calculations lead to  $\dot{V}_2 \le  -\frac{2\rho}{\mu_{\max}}V_2$,
 which further indicates 
$$
\lVert \psi \left( t \right) \rVert ^2\le \frac{\mu _{\max}}{\mu _{\min}}\lVert \psi \left( 0 \right) \rVert ^2e^{\frac{-2\rho}{\mu _{\max}}}.
$$
In summary, combined with Lemma \ref{lem2},  we can see that $x$ converges to $x^*$,  which is the optimal solution of problem (\ref{op2}),  with the exponential rate no less than $\frac{\rho}{\mu_{\max}}.$\hfill $\Box$

\begin{remark}
	From the perspective of convergence, algorithm (\ref{al2}) is able to ensure global exponential convergence of the system, while algorithm (\ref{al1})  can only achieve asymptotic convergence; therefore, algorithm (\ref{al2})  is superior to algorithm (\ref{al1}). However, the cost function has to be differentiable in algorithm (\ref{al2})  and the case with local constraints can not be coped with, while algorithm (\ref{al1}) is applicable to the non-smooth resource allocation problem with local feasibility set constraints. Additionally, the parameters in algorithm (\ref{al1}) have more selection margins than those in algorithm (\ref{al2}). Consequently, in terms of the application range of the algorithm, algorithm (\ref{al1}) is more widely applicable than algorithm (\ref{al2}).
\end{remark}

\section{Simulations}\label{sec.s}
In this section, for supporting the theoretical results, two examples are presented to verify and visualize the proposed algorithms.

\begin{figure}
	\centering
	\begin{tikzpicture}[->,>=stealth]
		\tikzstyle{node}=[shape=circle,fill=green!20,draw=none,text=black,inner sep=3pt, minimum size = 20 pt, ball color = green!50]
		\tikzstyle{input}=[
		shape=rectangle,
		fill=blue!20,
		draw=none,
		text=black,
		]
		
		\node[node] (1) at (0,0) {$1$};
		\node[node] (2) at (2.5,0) {$2$};
		\node[node] (3) at (5,0) {$3$};
		\node[node] (4) at (5,1.75) {$4$};
		\node[node] (5) at (2.5,1.75) {$5$};
		\node[node] (6) at (0,1.75) {$6$};
		
		\path (1) edge[->] (2)
		(2) edge[->] (3)
		(3) edge[->] (4)
		(4) edge[->] (5)
		(5) edge[->] (6)
		(6) edge[->] (1);
		
	\end{tikzpicture}
	\caption{Communication graph among  six generators.
	}\label{fig.g}
\end{figure}
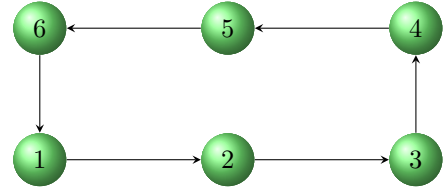

\textbf{Example 1:}
We apply algorithm (\ref{al1}) to the economic dispatch problem commonly found in smart grids, as done in \cite{Deng2020}. Consider a smart grid with six generators,  in which the communication digraph is described as Fig.\ref{fig.g}. The general economic dispatch problem requires determining the optimal output power for each generator to meet the network-wide load demand under security constraints.  Specifically,  in this instance, each generator is expected to cooperatively solve the constrained optimization problem as follows:
\begin{align}\label{exp1}
	&\min _{p_{G} \in \mathbb {R}^{6}} f(p_{G})=\sum _{i =1}^{N} f_{i} (p_{Gi}) \nonumber \\&\text {subject to} ~ \sum _{i =1}^{N} p_{Gi} = \sum _{i =1}^{N} p_{di} \\
	&\qquad  ~p_{Gi}^{\min } \le p_{Gi}  \le p_{Gi}^{\max }\nonumber
\end{align}
where $f_i  (p_{Gi})$ denotes the local generation cost function in $M\$$, $p_{di}$ is taken to represent the output power for generator $i$ in $MW$, and $p_{Gi}^{\max }$, $p_{Gi}^{\max }$ mean the upper and lower bound interpreted as the local bounded constraint of power generation. The local cost function $f_i(p_{Gi}) $ is defined as the following quadratic form:
\begin{align*}
	f_i(p_{Gi})=\gamma_i p_{Gi}^2  +\beta_i |p_{Gi} -\mathfrak{c}_i|+ \alpha_i,  \; \;  i\in \{1,\ldots,6 \}.
\end{align*}
It is obvious that $f_i(p_{Gi})$ is nonsmooth.

\begin{figure}
	\centering
	\includegraphics[width=9cm]{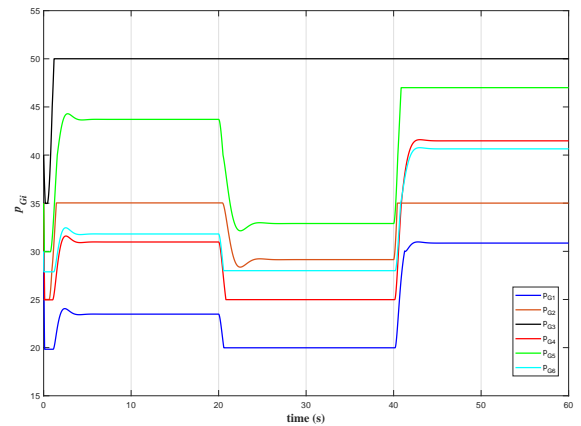}\\
	\caption{Evolutions of $p_{Gi}(t)$ in Example 1 with the load changed at $t=20s$ and $t=40s$.}\label{fig.1}
\end{figure}

Notice that problem (\ref{exp1}) is a particular case of problem (\ref{op1}). In problem (\ref{exp1}), the output power of each generator has to satisfy the coupled equality constraint and the local bounded constraint, while seeking the optimum distributively by sharing information with neighbors under the condition that only local data are available.
The involved parameters  $\alpha_i, \beta_i, \gamma_i, \mathfrak{c}_i, p_{Gi}^{\min }$  and $p_{Gi}^{\min }$ are shown in TABLE \ref{table:parameters}.
\begin{table}
	\begin{center}
		\caption{Parameters setting  \cite{Deng2020}}\label{table:parameters}
		\begin{tabular}{ccccccc}
			\hline
			$Generator$& $1$ & $2$ & $3$ & $4$ & $5$ & $6$\\ \hline
			$\alpha_i$ & $0.5$ & $1.5$    & $3$ & $1$ & $2.5$  & $2$\\
			$\beta_i$ & $3$ & $4$    & $5$ & $2$ & $3.5$ & $4.5$\\
			$\gamma_i$ & $2$ & $1$    & $0.5$ & $1.5$ & $1$ & $1.5$\\
			$\mathfrak{c}_i$ & $30$ & $28$   & $45$ & $35$ & $40$ & $35$\\
			$P_{di} $ & $45$ & $40$   & $25$ & $35$ & $30$ & $40$ \\
			$p_{Gi}^{\min }$ & $20 $ & $25 $   & $35$ & $25$ & $30 $ & $28$\\
			$p_{Gi}^{\max }$ & $40$ & $35$   & $50$ & $45$ & $47$ & $42$\\
			\hline
		\end{tabular}
	\end{center}
\end{table}
To run algorithm (\ref{al1}), we select the initial powers $p_{G}(0)=(30, 25, 40, 35, 35, 30)^T$, choose the auxiliary variables $\eta_i(0)=0, i \in \{1,\ldots,4\}$ and set the parameters $k_1=9, k_2=326$ and $k_3=5$. Since the local load demand $p_{di}$ may change in the actual application, we make $p_{d6}$ change from 40 $MW$ to 10 $MW$ at the time $t=20s$ in the simulation, and $p_{d6}$ change from 10 $MW$ to 70 $MW$ at $t=40s$.
The corresponding simulation results are displayed in Fig.\ref{fig.1}.

\begin{figure}
	\centering
	\includegraphics[width=9cm]{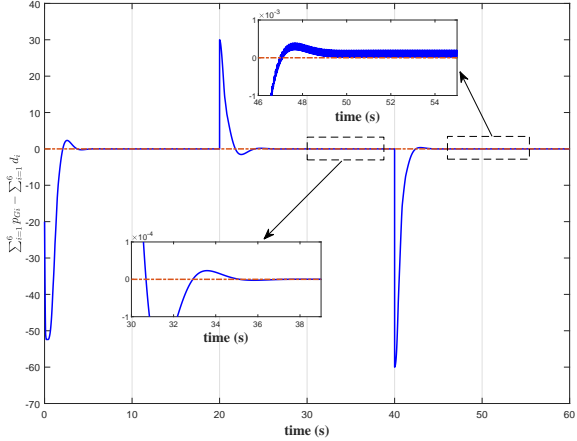}\\
	\caption{Evolutions of total mismatch in Example 1.}\label{fig.2}
\end{figure}

From Fig.\ref{fig.1},  it can be seen that the output powers of generators are always within the constraint sets, signifying that the local constraints are satisfied. In addition, the total mismatch is illustrated in Fig.\ref{fig.2}, and we can observe that it converges to 0 with high accuracy, indicating that the final output power of the generator satisfies the network resource constraint.

\begin{figure}
	\centering
	\includegraphics[width=9cm]{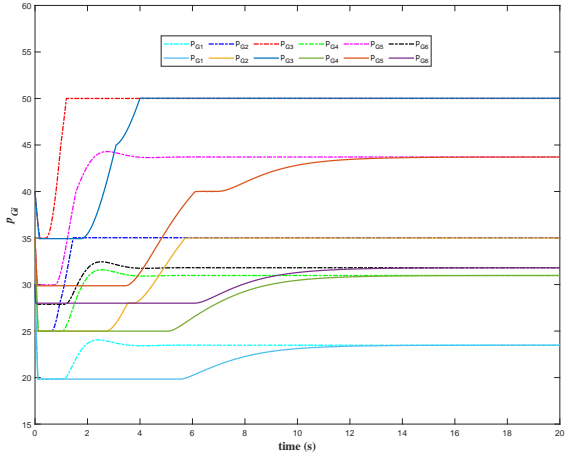}\\
	\caption{Comparison of simulation results between the algorithm in \cite{Deng2018} and algorithm (\ref{al1}) indicated by solid and dot dash lines, respectively.}\label{fig.3}
\end{figure}

Moreover, as a comparison to \cite{Deng2018}, for the same problem, we execute the algorithm in \cite{Deng2018}, setting the parameters as $k_1=2$ and $k_2=200$, and compare the simulation results with our proposed algorithm as depicted in Fig.\ref{fig.3},  it can be observed that both algorithms converge to the exact optimal solution of problem (\ref{exp1}), and our proposed algorithm has a faster convergence rate.

\textbf{Example 2:} Within this example, we verify the correctness of the convergence with regard to the algorithm described in Theorem \ref{the2} by means of more complex cost functions. Consider a network system with four agents (see \cite{DengHigh}), the communication graph among agents is a directed ring, and the cost functions are listed as follows: 
\begin{align*} {f_1}(x_1) =& x_{11}^2 + x_{12}^2 \\ {f_2}(x_2) =& \frac{x_{21}^2}{20 {x_{21}^2 + 1} } + \frac{x_{22}^2}{20 {x_{22}^2 + 1} } + ||x_2||^2 \\ {f_3}(x_3) =& ||x_3 - {[2\,3]^T}||^2 \\ {f_4}(x_4) =& \ln \left({{e^{ - 0.05{x_{41}}}} + {e^{0.05{x_{41}}}}} \right) \\ & + \ln \left({{e^{ - 0.05{x_{42}}}} + {e^{0.05{x_{42}}}}} \right) +||x_4||^2,
\end{align*}
where $x_i=(x_{i1}, \; x_{i2})  \in \mathbb{R}^2$.  Besides, the local resources of agents are denoted by $d_1=(2, \; 1)^T$, $d_2=(2, \; 3)^T$, $d_3=(2, \; 4)^T$, and $d_3=(1, \; 5)^T$,  respectively.
It is readily deduced that the mentioned cost functions are all differentiable, strongly convex, and the gradients of them are Lipschitz continuous, thus satisfying the assumptions of Theorem \ref{the2}.

\begin{figure}
	\centering
	\includegraphics[width=9cm]{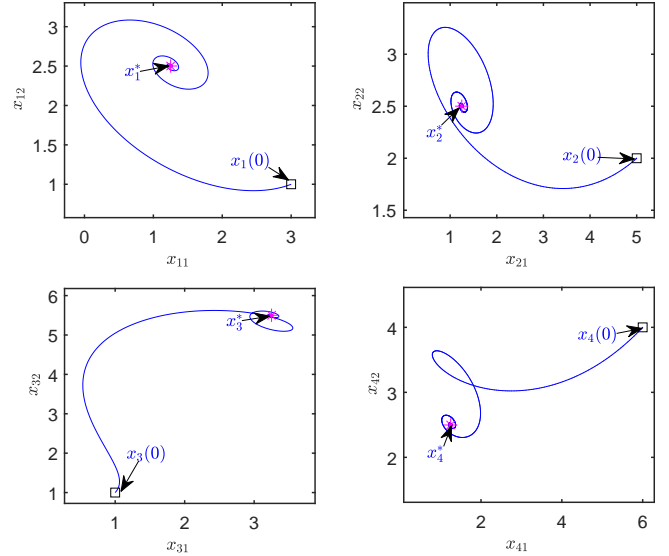}\\
	\caption{Evolutions of the states of four agents in Example 2.}\label{fig.4}
\end{figure}

To execute algorithm (\ref{al2}), we set parameters as $k_1=17,  k_2=290,  k_3=10$, and choose the auxiliary variables $\eta_i(0)=(0,\, 0)^T, i \in \{1,\ldots,4\}$ take the initial points as $x_1(0)=(3, \, 1)^T,$ $x_2(0)=(5, \, 2)^T,$  $x_3(0)=(1, \, 1)^T,$ $x_4(0)=(6, \, 4)^T.$ The simulation result is depicted in Fig.\ref{fig.4}, where the initial state is indicated by ``$\Box$" and ``$*$" indicates the final state of the agent.

\begin{figure}
	\centering
	\includegraphics[width=9cm]{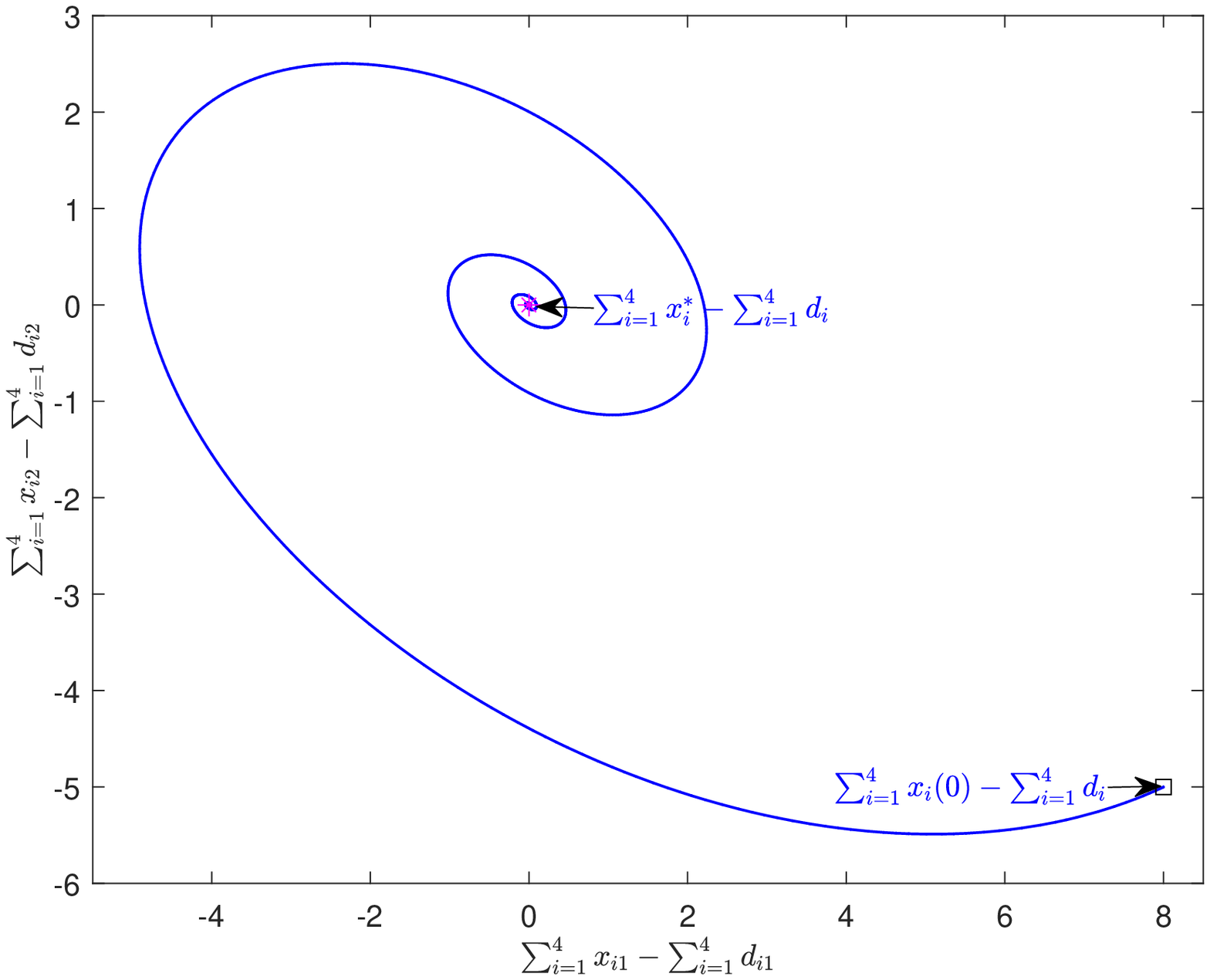}\\
	\caption{Evolutions of total mismath in Example 2.}\label{fig.5}
\end{figure}

It can be observed that the agent can eventually converge to the optimal solution. Fig.\ref{fig.5} displays the evolution trajectory curve of total mismatch, which can be seen to finally converge to the point $(0,\, 0)$, that is, the final decision satisfies the network resource constraint. In summary, the simulation results validate the effectiveness of algorithm (\ref{al2}).

\section{Conclusion}\label{sec.c}
This paper investigates the distributed nonsmooth resource allocation problem on multi-agent networks with heterogeneous set constraints that allow agents to communicate over a strongly connected and weight-balanced digraph. 
A continuous-time distributed algorithm is designed based on differential inclusion and differentiated projection operators. 
Further, the asymptotic convergence of the algorithm is proved by the set-valued LaSalle invariance principle and nonsmooth analysis theory, which allows the agent to achieve a globally optimal allocation while only exchanging local information with its neighbors. 
Moreover, the algorithm can achieve exponential convergence to the exact optimal allocation without involving local set constraints and with Lipschitz gradients. In order to verify the validity of algorithms, some simulation results are given at the end.
 In future work, further attention may be given to the design of algorithms on weight-unbalanced digraphs and the extension of algorithms to second-order or high-order systems.

\ifCLASSOPTIONcaptionsoff
  \newpage
\fi

\bibliographystyle{IEEEtran}
\bibliography{reso}

%
%
%
%

\end{document}